\documentclass[12pt]{amsart}

\usepackage[latin1]{inputenc}
\usepackage[T1]{fontenc}
\usepackage{hyperref}
\usepackage{geometry} 
\usepackage{stmaryrd}
\usepackage{amsmath,amssymb,amsthm,latexsym}

\newtheorem{theorem}{Theorem}%[section]
\newtheorem{proposition}[theorem]{Proposition}
\newtheorem{fact}[theorem]{Fact}
\newtheorem{cor}[theorem]{Corollary}
\newtheorem{lmm}[theorem]{Lemma}
\newtheorem*{claim}{Claim}
\newtheorem{clm}{Claim}

\theoremstyle{definition}
\newtheorem{definition}[theorem]{Definition}
\newtheorem{remark}[theorem]{Remark}
\newtheorem{conj}[theorem]{Conjecture}
\newtheorem{expl}[theorem]{Example}

\def\bsp{\begin{expl}}
\def\ebsp{\end{expl}}
\def\behe{\begin{clm}}
\def\ebehe{\end{clm}}
\def\beh{\begin{claim}}
\def\ebeh{\end{claim}}
\def\defn{\begin{definition}}
\def\edefn{\end{definition}}
\def\satz{\begin{theorem}}
\def\esatz{\end{theorem}}
\def\tats{\begin{fact}}
\def\etats{\end{fact}}
\def\kor{\begin{cor}}
\def\ekor{\end{cor}}
\def\bema{\begin{remark}}
\def\ebema{\end{remark}}
\def\lem{\begin{lmm}}
\def\elem{\end{lmm}}
\def\bem{\begin{remark}}
\def\ebem{\end{remark}}
\def\verm{\begin{conj}}
\def\everm{\end{conj}}
\def\bew{\begin{proof}}
\def\ebew{\end{proof}}
\def\bewbeh{\begin{proof}[Proof of Claim]}
\def\satzli{\begin{proposition}}
\def\esatzli{\end{proposition}}

\def\Ind#1#2{#1\setbox0=\hbox{$#1x$}\kern\wd0\hbox to 
0pt{\hss$#1\mid$\hss}
\lower.9\ht0\hbox to 0pt{\hss$#1\smile$\hss}\kern\wd0}
\def\ind{\mathop{\mathpalette\Ind{}}}
\def\Notind#1#2{#1\setbox0=\hbox{$#1x$}\kern\wd0\hbox to 0pt{\mathchardef
\nn="3236\hss$#1\nn$\kern1.4\wd0\hss}\hbox to 0pt{\hss$#1\mid$\hss}\lower.9\ht0
\hbox to 0pt{\hss$#1\smile$\hss}\kern\wd0}

\def\indd{\mathop{\ \hbox to 
0pt{$\mid^{\text{d}}$\hss}\,
\lower4pt\hbox to 0pt{\hss$\smile$\hss}\ \ }}
\def\nindd{\mathop{\ \hbox to 
0pt{$\!\not{\mid}^{\text{\,d}}$\hss}\,
\lower4pt\hbox to 0pt{\hss$\smile$\hss}\ 
\ }}

\def\M{\mathfrak M}

\def\ZZ{\mathbb Z}

\def\tp{\mathrm{tp}}

\def\ale{\lesssim}

\def\stab{\mathrm{Stab}}
\def\Sol{\mathfrak{S}}
\def\Zen{\mathfrak{Z}}

\begin{document}

\title{dp-minimal groups}          
\author{Frank O. Wagner}
\date{13 March 2024}
\address{Universit\'e Lyon 1; CNRS; Institut Camille Jordan UMR 5208, 21
avenue
Claude Bernard, 69622 Villeurbanne-cedex, France}
\email{wagner@math.univ-lyon1.fr}

\begin{abstract} A dp-minimal group is virtually nilpotent.
\end{abstract}
\thanks{Partially supported by ANR AAPG2019 GeoMod}
%\date{16/2/2015}

\maketitle

\defn A partial type $\pi(x)$ is dp-minimal if there are no formulas $\varphi(x,y)$ and $\psi(x,y)$ and mutually indiscernible sequences $(a_i)_{i<\omega}$ and $(b_i)_{i<\omega}$ such that for any $i,j\in\omega$ the type 
$$\pi(x)\cup\{\varphi(x,a_i),\psi(x,b_j)\}\cup\{\neg\varphi(x,a_s):s\not=i\}\cup\{\neg\psi(x,b_t):t\not=j\}$$ is consistent.

Equivalently, for any two mutually indiscernible sequences $(a_i)_{i<\omega}$ and $(b_i)_{i<\omega}$ and any $c\models\pi$, either $(a_i)_{i<\omega}$ or $(b_i)_{i<\omega}$ is indiscernible over $c$.\edefn

Recall that a dp-minimal theory is dependent, so in particular the Baldwin-Saxl condition holds, and any intersection of uniformly definable subgroups is equal to a (uniformly) finite subintersection. Moreover, Pierre Simon has shown \cite[Claim in Proof of 4.31]{Si15} that in a dp-minimal group, for any two definable subgroups $H$ and $K$ the intersection $H\cap K$ has finite index in either $H$ or in $K$. By compactness, for any family of uniformly definable subgroups $\{H_i:i\in I\}$ there is $n<\omega$ such that $H_i\cap H_j$ has index at most $n$ in either $H_i$ or in $H_j$, for any $i,j\in I$.

We shall call a group {\em $\bigwedge$-definable} if it is an intersection of definable groups. Recall that a group us {\em virtually P} if it has a subgroup of finite index which is P.

\bem If $G$ is virtually soluble, then its soluble radical $R(G)$ (the subgroup generated by all normal soluble subgroups) is soluble and definable; if $G$ is virtually nilpotent, its Fitting subgroup $F(G)$ (the subgroup generated by all normal nilpotent subgroups) is nilpotent and definable. If $G$ is soluble, it has a definable characteristic series with abelian quotients; if $G$ is nilpotent, it has a definable characteristic series with central quotients.\ebem
\bew If $G$ has a soluble subgroup $H$ of finite index, $H$ contains a normal soluble subgroup $N$ of finite index. Now for any normal soluble subgroup $S$ the group $NS$ is again normal soluble. Hence $G$ has a maximal normal soluble subgroup of finite index, so $R(G)$ is soluble. It is then definable by \cite[Theorem 1.1]{OH13}. The proof for $F(G)$ is analogous, replacing soluble by nilpotent.

Now if $G$ is soluble of derived length $n+1$, define $c_0(x)=x$ and $c_{i+1}(x_1,\ldots,x_{2^{i+1}})=[c_i(x_1,\ldots,x_{2^i}),c_i(x_{2^i+1},\ldots,x_{2^{i+1}})]$. Then 
$$A=Z(C_G(c_n(x_1,\ldots,x_{2^n}):x_1,\ldots,x_{2^n}\in G)$$
is a definable characteristic abelian subgroup containing $G^{(n)}$. We divide out by $A$ and finish by induction. For nilpotency, the claim is obvious, as $Z(G)$ is clearly definable.\ebew

\satzli A dp-minimal group $G$ is virtually soluble.\esatzli
\bew By a theorem of Simon \cite[Proposition 4.31]{Si15}, $G$ has a definable characteristic abelian subgroup $A$ such that $G/A$ has finite exponent.
As $G/A$ is still dp-minimal, we may assume that $G$ has finite exponent.

Let $H\le G$ be a definable subgroup. The family of its $G$-conjugates is a family of uniformly definable subgroups, so there is $n<\omega$ such that $H\cap H^g$ has index $\le n$ either in $H$ or in $H^g$. Suppose $|H:H\cap H^g|\le n$ (the other case is analogous). Then $|H^g:H^g\cap H^{g^2}|\le n$, and $|H\cap H^g:H\cap H^g\cap H^{g^2}|\le n$, whence $|H:H\cap H^g\cap H^{g^2}|\le n^2$. Iterating, we see that $|H:H\cap\cdots\cap H^{g^k}|\le n^k$. But for $k=o(g)-1$ the subgroup $N=H\cap\cdots\cap H^{g^k}$ is $g$-invariant, so $|H:N|=|H^g:N|\le n^k$. But this implies $|H:H\cap H^g|=|H^g:H\cap H^g|\le n$. Now NIP implies that there is $\ell$ such that the intersection of any finite number of $G$-conjugates of $H$ is an intersection of at most $\ell$ of them, and has index $\le \ell n$ in $H$. Hence $H^{nor}=\bigcap_{g\in G}H^g$ is a normal subgroup of index at most $\ell n$ in $H$, definable over the same parameters as $H$.

Let $\M$ be a model, and $\Sol$ the set of $M$-definable normal soluble subgroups. Since $\Sol$ is closed under product, $\bigcup\Sol$ is a normal subgroup of $G$.
If $G$ is not virtually soluble, there is a minimal $M$-$\bigwedge$-definable subgroup $H$ which is not virtually contained in any subgroup in $\Sol$. Let $p$ be an f-generic type of $H$ over $M$. Note that by minimality $H=H^0$, so $H$ is connected; if $H=\bigcap_i H_i$ then $H=\bigcap_i H_i^{nor}$ as well by connectivity, so $H$ is an intersection of $M$-definable normal subgroups and must itself be normal.

By dp-minimality, $p$ is either generically stable or distal. If $p$ is distal, then by the generalization of \cite[Section 5]{St23} to $\bigwedge$-definable normal subgroups of a dp-minimal group, $H$ is nilpotent. By compactness it is contained in a definable normal nilpotent group of the same class which must be in $\Sol$, a contradiction.

So $p$ is generically stable. Let $(a_i)_\ZZ\hat\ (b_i)_\ZZ$ be a Morley sequence in $p$ over $M$.

If $[a_0,b_0]$ lies in some $M$-definable coset of some group $S\in\Sol$, say $[a_0,b_0]\in mS$, then $b_0$, $b_1$ and $b_0b_1$ all satisfy $p|Ma_0$ by generic stability and genericity of $p$, so
$$[a_0,b_0]=[a_0,b_1]=[a_0,b_0b_1]=[a_0,b_1][a_0,b_0]^{b_1}\mod S.$$
Hence $[a_0,b_0]^{b_1}\in S$, so $[a_0,b_0]\in S$. Now $C_G(a_0/S)$ contains a generic $b_0$, so has finite index in $H$, which is connected. By compactness $C_G(a_0/S)$ contains an $M$-definable normal supergroup $K$ of $H$. Then $Z(K/S)$ is $M$-definable and contains $a_0$ which is generic in $H$. By connectivity $H\le Z(K/S)$. But $Z(K/S)\in\Sol$, a contradiction.

In particular $q=\tp([a_0,b_0]/M)$ is a non-algebraic type over $M$; it is generically stable since $\tp(a_0,b_0/M)$ is. By \cite[Lemma 6.2]{St23} the stabilizer $\stab(q)$ is infinite; in fact it contains $c_0c_1^{-1}$ for some/any independent $c_0,c_1\models q$. Thus the coset $\stab(p)c_1$ does not depend on $c_1\models q$ and must be $M$-definable; moreover it contains all realizations of $q$. By the previous paragraph, $\stab(q)\notin\Sol$; as $q(x)$ implies $x\in H$ we get $\stab(q)\le H$ and we must have equality by minimality of $H$.

Work over $M$. By dp-minimality of $q$, one of $(a_i)_\ZZ$ or $(b_i)_\ZZ$ is indiscernible over $[a_0,b_0]$, whence independent of $[a_0,b_0]$. In particular $[a_0,b_0]$ is independent of $a_0$ or of $b_0$; by symmetry $[b_0,a_0]$ is independent of $b_0$ or of $a_0$. But $[b_0,a_0]=[a_0,b_0]^{-1}$, so it is independent of $a_0$ and of $b_0$. It follows that
$$H=\stab(q)=\stab(q|Mb_0)=\stab([a_0,b_0]/Mb_0)=\stab(b_0^{-a_0}b_0/Mb_0)=\stab(b_0^{-a_0}/Mb_0).$$
But $H=\stab(p)=\stab(p|Mb_0)$ by generic stability of the generic type $p$. Then $H$ is the {\em right} stabilizer of $p^{-1}$. So for two independent realizations $a,a'\models p$ we have that 
$$\tp(a^{-1}a'/Ma')=p^{-1}a'|Ma'=p^{-1}|Ma'\quad\mbox{and}\quad\tp(a^{-1}a'/Ma)=a^{-1}p|Ma=p|Ma.$$
In particular $p=\tp(a^{-1}a'/M)=p^{-1}$, and $H$ is the left and right stabilizer of $p$. Similarly, if $a\models p|Mb_0$ and $c\models\tp(b_0^{-a_0}/Mb_0)$ with $a\ind_{Mb_0}c$, then 
$$p|Mb_0=\tp(ac/Mb_0)=\tp(b_0^{-a_0}/Mb_0).$$
Since $b_0^{-1}\models d_px\,\exists z\in G\,x=y^z$, there is $m'\in M$ such that $p\vdash x\in m'^G$. As $H=p\cdot p$ is normal in $G$, we have that $H=m'^G\cdot m'^G$ is $M$-definable. Again by definability of $p$ there is $m\in M$ such that $p\vdash x\in m^H$. (In fact it is easy to see that we can choose $m=m'$.)

It follows that $m\in H(M)\setminus\bigcup\Sol(M)$; in particular $H(M)/\bigcup\Sol(M)$ is nontrivial. Note that if there is $S\in\Sol$ such that $m^2\in S$, then every realization of $p$ has order two modulo $S$. But for two independent realizations $a,a'\models p$ we have $aa'\models p$, whence $a^2, a'^2, (aa')^2\in S$ and
$$aa'=(aa')^{-1}=a'^{-1}a^{-1}=a'a\mod S.$$
So $H/S$ is generically commutative and $[a_0,b_0]\in S$; as above $H\le\bigcup\Sol$, a contradiction.
It follows that $m^2\notin\bigcup\Sol$.

Now $m\not=m^{-1}\models p^{-1}=p$, so there is $h\in H$ with $m^{-1}=m^h$. Since $m\in M$ we can choose $h\in H(M)$. Clearly $h\notin\bigcup\Sol$, and $m\in C_H(h^2)\setminus C_H(h)$.

Consider any $c\in H(M)\setminus\bigcup\Sol(M)$, and put $q'=\tp(c^{a_0}/M)$, a generically stable type. By \cite[Lemma 6.2]{St23} for $T=\stab(q')$ and any $t\models q'$ the coset $Tt$ does not depend on $t$, is $M$-definable and contains $q'$.

If $T$ were soluble, then for generic $a\models p|Ma_0$ we have $\tp(a_0a/M)=p=\tp(a_0/M)$, so $q'^a=q'\subseteq (Tt)^a\cap Tt$, whence $[a,t^{-1}]\in T$ (recall that $T$ is a stabilizer, whence connected, and thus normal). Hence $C_H(Tt/T)$ is a generic $M$-definable subgroup of $H$, and must equal $H$ by connectivity. Thus $t$ is central in $H/T$, as is $t^G$, and there is a definable $S\in\Sol$ containing $Tt$. Then $c\in S$, a contradiction.

Since $T\le H$ we get $T=H$ by minimality. But then $q'=p$ as before, and $c\in m^H$.
So $H(M)/\bigcup\Sol(M)$ is a group of order at least $3$, finite exponent, and a single non-trivial conjugacy class, a contradiction.\ebew

\satzli Let $G$ be a soluble dp-minimal group. Then $G$ is virtually nilpotent.\esatzli
\bew By \cite[Theorem 5.21]{St23} we may assume that $G$ has a generically stable principal f-generic type. Since $G$ is soluble, it has a definable normal series with abelian quotients. So by induction on the derived length we may assume that there is a normal definable nilpotent subgroup $N$ such that $G/N$ is finite-by-abelian. Let $F/N$ be
the finite subgroup such that $G/F$ is abelian. Then $\bar n^g\in \bar nF$ for every $\bar n\in G/N$, and $|G:C_G(\bar n)|\le|F|$. Since $G$ is NIP, the intersection of all $C_G(\bar n)$ is a finite subintersection $G^c$. We replace $G$ by $G^c$ and assume that $G/N$ is abelian.

Let $\M$ be an $\omega$-saturated model. Recall that $G$ has a characteristic $\emptyset$-definable abelian subgroup $A$ such that $G/A$ has finite exponent. Put 
$$\Zen=\{Z\le G:Z=Z_n(G_Z)\mbox{ for some normal $M$-definable $G_Z\le G$ of finite index}\}.$$ Note that $Z_n(G_Z)\,Z_m(G_{Z'})\le Z_{n+m}((G_Z\cap G_{Z'})ZZ')$, so $\Zen$ is closed under product, and consists of normal $\emptyset$-definable groups. We want to show that $N\le Z$ for some $Z\in\Zen$. So suppose not.

Let $H$ be a minimal $M$-$\bigwedge$-definable subgroup not virtually contained in $\bigcup\Zen$. If $A$ is not virtually contained in $\bigcup\Zen$ we choose $H\le A$. 
%Note that if there were $Z=Z_n(G_Z)\in\Zen$ such that $H/Z$ is finite, then $H\le Z_{n+1}(G_0)$, where $G_0=C_G{G_Z}(H/Z)$ has finite index in $G$, a contradiction. 
Note that $H$ is connected, and $G$-invariant, as $G/C_G(H)$ has finite exponent, say $k$. Dividing out by some $Z\in Z$ we may assume that $H$ is abelian.

Consider a principal generic element $g$ of $G$ over $M$, and a $(g,M)$-definable subgroup $K_g$ of $G$. Then for an independent $g'\equiv_M g$ we either have $M_g\ale M_{g'}$ or $M_{g'}\ale M_g$. But by generic stability of $\tp(g/M)$ we have $g,g'\equiv_M g',g$, so $M_g$ and $M_{g'}$ are (uniformly) commensurable. For any $g''\equiv_M g$ we choose $g'\equiv_M g$ independent of $g,g''$ over $M$, so $M_g$ and $M_{g'}$ are both uniformly commensurable to $M_{g''}$, and hence uniformly commensurable. It follows that the intersection $\bigcap_{g'\equiv_M g}M_{g'}$ equals a finite subintersection, has finite index in $M_g$ and is $M$-definable.

In particular, for any $Z\in\Zen$ the centraliser $C_G(g/Z)$ has an $M$-definable subgroup $K$ of finite index normal in $G$ which centralises modulo $Z$ some normal generic subgroup $G_K\le G$. Then $K\cap G_K$ has finite index in $C_G(g/Z)$, and is contained in $\Zen$. Moreover, the subgroup $[g,H]$ is $(g,M)$-$\bigwedge$-definable (as $H$ is abelian) and connected (since $H$ is), whence $M$-$\bigwedge$-definable; it must either equal $H$ or be virtually contained in $\bigcup\Zen$, whence contained in $\bigcup\Zen$ by connectedness. The latter case implies that $H\le\bigcup\Zen$ again by connectedness, a contradiction.

Let $R$ be the ring of quasi-endomorphisms modulo equivalence of $H/\bigcup\Zen$ generated by $G$, and their inverses (where they exist), in the monster model. Note that by the above, the map $r_g:x\mapsto [g,x]$ has kernel virtually contained in $\bigcup\Zen$ and is surjective, whence invertible in $R$. Composing with conjugation by an independent principal generic $g'$, we see that $f_{g,g'}:x\mapsto x^{-gg'}x^{g'}$ is also invertible in $R$.

Now let $I$ be a maximal ideal in $R$, so $R/I$ is a field. Since conjugation by $g$ has order dividing $k$, there are only $k$ different images in $R/I$ for conjugation by realizations of $p$.
But since $f_{g,g'}$ is invertible, it cannot be in $I$. As there is a unique type of two independent realizations of $p$, which is realized by $gg'$ and $g$, independent generics have a different image in $R/I$, a contradiction.

This finishes the proof.
\ebew


\begin{thebibliography}{99}
\bibitem{BH79} Bryant, R. M. and B. Hartley, {\em Periodic locally soluble groups with the minimal condition on centralizers}, J. Alg.\ 61:328--334, 1979.
\bibitem{Ke89} Kegel, O., {\em Four lectures on Sylow theory in locally finite groups}, pp. 3-28 in {\em Group Theory}, edited by K. N. Nah and Y. K. Leong, Walter de Gruyter, Amsterdam, 1989.
\bibitem{OH13} Ould Houcine, A., {\em A remark on the definability of the Fitting subgroup and the soluble radical}, Math.\ Log.\ Quart.\ 1--4, 2013.
\bibitem{Si15} Pierre Simon. {\em A Guide to NIP Theories}. ASL Lecture notes in Logic 44, CUP 2015.
\bibitem{St23} Atticus Stonestrom. {\em On non-abelian dp-minimal groups}, arXiv:2308.04209

\end{thebibliography}
\end{document}